\documentclass[12pt]{amsart}
\advance\hoffset by -0.5 truecm
\usepackage{amssymb}
\newtheorem{Theorem}{Theorem}[section]

\newtheorem{Proposition}[Theorem]{Proposition}

\newtheorem{Remark}[Theorem]{Remark}

\def\V{\mbox{Var}}

\def\Z{{\mathbb Z}}
\def\R\re
\def\V{\bf V}

\def \la{\lambda}

\def \re{{\mathbb R}}

\def \0{\lambda_{0}}
\def \la{\lambda}

\begin{document}
\title[KAM-cocycle of a magnetic flow]{The longitudinal KAM-cocycle of a magnetic flow}

\author[G. P. Paternain]{Gabriel P. Paternain}
%\thanks{G. P. Paternain
% was partially
% supported by CIMAT, Guanajuato, M\'{e}xico}
 \address{ Department of Pure Mathematics and Mathematical Statistics,
University of Cambridge,
Cambridge CB3 0WB, England}
 \email {g.p.paternain@dpmms.cam.ac.uk}

%\author[J. Petean]{Jimmy Petean}
% \address{CIMAT  \\
%          A.P. 402, 36000 \\
%          Guanajuato. Gto. \\
%          M\'{e}xico.}
%\email{jimmy@cimat.mx}

%\thanks{J. Petean is supported by grant 37558-E of CONACYT}

%\subjclass{53C25, 53C21, 58F17, 35J15}

\date{May 2004}

%\maketitle

\begin{abstract} Let $M$ be a closed oriented surface of negative Gaussian curvature and let
$\Omega$ be a non-exact 2-form. Let $\lambda$ be a small positive
real number. We show that the longitudinal KAM-cocycle of the
magnetic flow given by $\la\,\Omega$ is a coboundary if and only
if the Gaussian curvature is constant and $\Omega$ is a constant
multiple of the area form.

\end{abstract}

\maketitle

\section{Introduction}

Let $M$ be a closed oriented surface with negative Gaussian
curvature. It is a classical fact that the geodesic flow $\phi$ of
$M$ is a contact Anosov flow (cf. \cite{KH,P1}). Being contact means that there
exists a contact form $\alpha$ on the unit sphere bundle $SM$ such
that the vector field $X$ on $SM$ that generates $\phi$ is
determined by the equations $\alpha(X)\equiv 1$ and
$i_{X}d\alpha\equiv 0$. The Anosov property means that $T(SM)$
splits as $T(SM)=E^{0}\oplus E^{u}\oplus E^{s}$ in such a way that
there are constants $C>0$ and $0<\rho<1<\eta$ such that $E^{0}$ is
spanned by $X$ and for all $t>0$ we have
\[\|d\phi_{-t}|_{E^{u}}\|\leq C\,\eta^{-t}\;\;\;\;\mbox{\rm
and}\;\;\;\|d\phi_{t}|_{E^{s}}\|\leq C\,\rho^{t}.\]

The subbundles are then invariant and H\"older continuous and have
smooth integral manifolds, the stable and unstable manifolds,
which define a continuous foliation with smooth leaves. The Anosov
property immediately implies that $E^{u}$ and $E^{s}$ must be
contained in the kernel of the contact form $\alpha$ and thus
$\mbox{\rm Ker}\,\alpha=E^{u}\oplus E^{s}$ and $E^{u}\oplus E^{s}$
is a $C^{\infty}$ subbundle.

Suppose now that $\la$ is a small positive number. Geodesics are
smooth curves in $M$ with zero geodesic curvature. Consider
instead smooth curves in $M$ with constant geodesic curvature
$\la$. We will call such curves {\it magnetic geodesics}. Through
each point $x\in M$, there is a unique magnetic geodesic with
velocity $v\in T_{x}M$, $|v|=1$. Thus magnetic geodesics define
also a flow $\phi^{\la}$ on $SM$ that we will call the {\it magnetic
flow} of $M$. For $\la$ small, $\phi^{\la}$ will still be Anosov,
although in general it will not be contact anymore. The flow
$\phi^{\la}$ is nevertheless still a Hamiltonian flow. A
straightforward calculation shows that $\phi^{\la}$ is the
Hamiltonian flow of the Hamiltonian
$H(x,v)=\frac{1}{2}|v|^{2}_{x}$ with respect to the symplectic
form on $TM$ given by
\[-d\alpha+\la\,\pi^{*}\Omega_{a},\]
where $\Omega_{a}$ is the area form of $M$ and $\pi:TM\to M$ is
the canonical projection. Hence, for small values of $\la$,
$\phi^{\la}$ is a volume preserving Anosov flow and in fact, it
preserves the volume form $\alpha\wedge d\alpha$.

To any $C^{k}$ volume preserving Anosov flow $\varphi$ on a closed
3-manifold $N$, P. Foulon and B. Hasselblatt \cite{FH} associated
its {\it longitudinal KAM-cocycle}. This is a cocycle that
measures the regularity of the subbundle $E^{u}\oplus E^{s}$ and
whose definition we now recall. Consider local coordinates (cf.
\cite{KH,FH}) adapted to the stable and unstable foliations
$\psi_{p}:(-\varepsilon, \varepsilon)^{3}\to N$, $p\in N$, and
denote the coordinate variables $(u,t,s)$. Consider the
transversals $\Delta_{p}:=\psi_{p}((-\varepsilon,
\varepsilon)\times\{0\}\times (-\varepsilon, \varepsilon))$ and
$\Delta_{\varphi_{T}(p)}:=\psi_{\varphi_{T}(p)}((-\varepsilon,
\varepsilon)\times\{0\}\times (-\varepsilon, \varepsilon))$. Then
$\Delta_{\varphi_{T}(p)}\cap \varphi_{T}(\Delta_{p})$ contains
local strong unstable and stable manifolds of $\varphi_{T}(p)$,
but the two transversals are not in general identical. Let
$f_{T}(u,s)$ be the time lengths of the orbit segments between
$\Delta_{\varphi_{T}(p)}$ and $\varphi_{T}(\Delta_{p})$ and set
\[K(p,T):=\frac{\partial^{2}f_{T}}{\partial u\partial s}(0,0).\]
Foulon and Hasselblatt show that $K(p,T)$ is an additive cocycle
(i.e. $K(p,T+S)=K(\varphi_{T}(p),S)+K(p,T)$) and that the
cohomology class of $K$ is independent of the adapted coordinates.
The main theorem in \cite{FH} asserts that $E^{u}\oplus E^{s}$ is 
always Zygmund-regular and that the following are
equivalent:
\begin{enumerate}
\item $E^{u}\oplus E^{s}$ is ``little Zygmund'';
\item the longitudinal KAM-cocycle is a coboundary;
\item $E^{u}\oplus E^{s}$ is Lipschitz;
\item $E^{u}\oplus E^{s}$ is $C^{k-1}$;
\item $\varphi$ is a suspension or contact flow.
\end{enumerate}

(A continuous function $f:U\to\re$ on an open set $U\subset \re$ is said to be
{\it Zygmund-regular}
 $|f(x+h)+f(x-h)-2f(x)|=O(h)$ for all $x$ in $U$. 
The function is said to be {\it ``little Zygmund''} if $|f(x+h)+f(x-h)-2f(x)|=o(h)$.)

It is well known that for flows, a ``choice of time" or equivalently, a choice of speed at which
orbits travel gets reflected on the regularity of the corresponding strong stable and strong unstable distributions. The situation is different if we look at the weak unstable and stable bundles $E^{0}\oplus E^{u}$ and
$E^{0}\oplus E^{s}$. S. Hurder and A. Katok proved \cite{HK} that
the weak bundles are always differentiable with Zygmund-regular derivative
and there is a cocycle obstruction to higher regularity given by the first
nonlinear term in the Moser normal form (this explains why Foulon and 
Hasselblatt used the terminology ``longitudinal KAM-cocycle'').

The question we would like to address here is: when is the
longitudinal KAM-cocycle of the magnetic flow $\phi^{\la}$ a
coboundary?

\medskip

\noindent {\bf Theorem A.} {\it Let $K$ be the Gaussian curvature
of $M$ and suppose $2\la^{2}+K(x)<0$ for all $x\in M$. Then the
longitudinal KAM-cocycle of $\phi^{\la}$, $\la\neq 0$, is a
coboundary if and only if $K$ is constant.

}

\medskip

We note that if $\la^{2}+K<0$, then $\phi^{\la}$ is Anosov (as it can be easily seen from the corresponding Jacobi equation) but there can be other larger values of $\la$ for which
$\phi^{\la}$ is Anosov (cf. \cite{BP}). Most likely Theorem A is
also true for any value of $\la$ for which $\phi^{\la}$ is Anosov,
but our methods do not yield that much.

There are earlier versions of Theorem A in the literature. Suppose
that we replace $\Omega_{a}$ by an arbitrary two-form $\Omega$ and
consider as above the Hamiltonian flow of $\frac{1}{2}|v|^{2}_{x}$
with respect to the symplectic form given by
$-d\alpha+\la\,\pi^{*}\Omega$. We proved in \cite{P} using
Aubry-Mather theory (and in any dimensions) that if $\Omega$ is
exact, then the corresponding magnetic flow has a $C^{1}$ Anosov
splitting only if $\la=0$, i.e. the flow is geodesic. Unfortunately, these
methods do not carry over to the non-exact case. In an unpublished
manuscript \cite{B}, J. Boland proved Theorem A under an
additional assumption on the metric. He showed that if there is a
closed geodesic transverse to which the Gaussian curvature is
infinitesimally increasing, then the Anosov splitting is not
differentiable for small $\la$, unless $\la=0$ (how small $\la$
had to be was unspecified). His methods were quite different from
ours and were based on a probabilistic approach using the
Feynman-Kac formula.

 Our proof will be based on
establishing results for magnetic flows analogous to Theorem 3.6 in \cite{GK}.
With little extra effort they will yield a more
general theorem as we now explain.

Suppose $\Omega$ is an arbitrary non-exact 2-form and $\phi^{\la}$
its associated magnetic flow.

\medskip

\noindent {\bf Theorem B.} {\it Let $M$ be a closed oriented
surface with negative Gaussian curvature. There exists $\la_{0}>0$
such that the longitudinal KAM-cocycle of $\phi^{\la}$ for
$0<|\la|<\la_{0}$ is a coboundary if and only if $K$ is constant
and $\Omega$ is a constant multiple of the area form.

}

\medskip

\section{Geometry of $SM$}

Let $M$ be a closed oriented surface, $SM$ the unit sphere bundle
and $\pi:SM\to M$ the canonical projection. The latter is in fact
a principal $S^{1}$-fibration and we let $V$ be the infinitesimal
generator of the action of $S^1$.

Given a unit vector $v\in T_{x}M$, we will denote by $iv$ the
unique unit vector orthogonal to $v$ such that $\{v,iv\}$ is an
oriented basis of $T_{x}M$. There are two basic 1-forms $\alpha$
and $\beta$ on $SM$ which are defined by the formulas:
\[\alpha_{(x,v)}(\xi):=\langle d_{(x,v)}\pi(\xi),v\rangle;\]
\[\beta_{(x,v)}(\xi):=\langle d_{(x,v)}\pi(\xi),iv\rangle.\]
The form $\alpha$ is precisely the contact form that we mentioned
in the introduction.

A basic theorem in 2-dimensional Riemannian geometry asserts that
there exists a unique 1-form $\psi$ on $SM$ (the connection form)
such that

\begin{align*}
\psi(V)&=1\\ & d\alpha=\psi\wedge \beta\\ & d\beta=-\psi\wedge
\alpha\\ & d\psi=-(K\circ\pi)\,\alpha\wedge\beta
\end{align*}
where $K$ is the Gaussian curvature of $M$. In fact, the form
$\psi$ is given by
\[\psi_{(x,v)}(\xi)=\left\langle \frac{DZ}{dt}(0),iv\right\rangle,\]
where $Z:(-\varepsilon,\epsilon)\to SM$ is any curve with
$Z(0)=(x,v)$ and $\dot{Z}(0)=\xi$ and $\frac{DZ}{dt}$ is the
covariant derivative of $Z$ along the curve $\pi\circ Z$.

It is easy to check that $\alpha\wedge\beta=\pi^{*}\Omega_{a}$,
hence
\begin{equation}
d\psi=-\pi^{*}(K\,\Omega_{a}). \label{psi}
\end{equation}

For later use it is convenient to introduce the vector field $H$
uniquely defined by the conditions $\beta(H)=1$ and
$\alpha(H)=\psi(H)=0$. The vector fields $X,H$ and $V$ are dual to
$\alpha,\beta$ and $\psi$.

\section{Proof of Theorems A and B} Let $\Omega$ be an arbitrary smooth
2-form. We write $\Omega=F\,\Omega_{a}$, where $F:M\to\re$ is a
smooth function.

 Since $H^{2}(M,\re)=\re$, there exist a constant $c$ and a
smooth 1-form $\theta$ such that
\[\Omega=cK\,\Omega_{a}+d\theta\]
and $c=0$ if and only if $\Omega$ is exact. Using (\ref{psi}) we
have
\[\omega_{\la}:=-d\alpha+\la\,\pi^{*}\Omega=d(-\alpha-\la
c\,\psi+\la\,\pi^{*}\theta).\] The vector field $X_{\la}$ that
generates $\phi^{\la}$ is given by $X_{\la}=X+\la F\,V$ since it
satisfies the equation $dH=i_{X_{\la}}\omega_{\la}$. If we
evaluate the primitive $-\alpha-\la c\,\psi+\la\,\pi^{*}\theta$ of
the symplectic form $\omega_{\la}$ on $X_{\la}$ we obtain:
\begin{equation}
(-\alpha-\la c\,\psi+\la\,\pi^{*}\theta)(X_{\la})(x,v)=-1-\la^{2}
F(x) c+\la\,\theta_{x}(v).\label{contact}
\end{equation}
Let us prove the easy part of Theorems A and B. Suppose $M$ has
constant curvature and $\Omega$ is a constant multiple of the area
form, i.e. $F$ is constant. Then we can choose $\theta=0$ and thus
$-\alpha-\la c\,\psi$ is a primitive of $\omega_{\la}$ which on
the vector field $X_{\la}$ is a constant equal to $-1-\la^{2} F
c$. It follows that $\phi^{\la}$ is a contact flow for all values
of $\lambda$ except those for which $-1=\la^{2} F c$ (in which
case the flow is in fact the horocycle flow). If the flow is
contact, then of course, its longitudinal KAM-cocyle is a
coboundary.

Suppose now that the longitudinal KAM-cocyle is a coboundary. By
the main theorem of Foulon and Hasselblatt that we mentioned in
the introduction, there is a smooth $\phi^{\la}$-invariant 1-form
$\tau$ which is equal to 1 on $X_{\la}$ and whose kernel is
$E^{u}\oplus E^{s}$. By ergodicity, there exists a constant $k$
such that
\[d\tau=k\omega_{\la}.\]
Hence
\[d(\tau+k\alpha+\la kc\,\psi-k\la\,\pi^{*}\theta)=0.\]
Since $\pi^{*}:H^{1}(M,\re)\to H^{1}(SM,\re)$ is an isomorphism,
there exists a closed 1-form $\delta$ on $M$ and a smooth function
$g:SM\to\re$ such that
\[\tau+k\alpha+\la
kc\,\psi-k\la\,\pi^{*}\theta=\pi^{*}\delta+dg.\] Evaluating both
sides on $X_{\la}$ we obtain
\begin{equation}
1+k+\la^{2} F(x)
kc-k\la\,\theta_{x}(v)=\delta_{x}(v)+dg(X_{\la}).\label{prekey}
\end{equation}
Since $X$ and $V$ preserve the volume form $\alpha\wedge d\alpha$,
then so does $X_{\la}=X+\la\,F V$ and thus $\phi^{\la}$ preserves
the normalized Liouville measure $\mu$ of $SM$. If we integrate
(\ref{prekey}) with respect to $\mu$ we get $$1+k+\la^{2}kc\int
F\,d\mu =0$$ since $\mu$ is invariant under the flip $v\mapsto
-v$. It follows that $k$ is a non-zero constant and when $\Omega$
is the area form (Theorem A) we are left with the equation
\begin{equation}
-k\la\,\theta_{x}(v)-\delta_{x}(v)=dg(X_{\la})=X_{\la}(g).\label{llave}
\end{equation}

Theorem A will now be a consequence of the following result which
we prove in the next section.

\begin{Theorem}Let $\Omega$ be the area form and suppose $2\la^{2}+K(x)<0$ for all $x\in M$.
If $\omega$ is any smooth 1-form on $M$ such that there is a
smooth function $g:SM\to\re$ for which
$\omega_{x}(v)=X_{\la}(g)$, then $\omega$ is exact.
\label{1forma}
\end{Theorem}

If we now apply Theorem \ref{1forma} to the form
$-k\la\,\theta-\delta$ we conclude that $\theta$ must be a closed
form. But if $\theta$ is closed, $\Omega_{a}=cK\,\Omega_{a}$ and
thus $K$ is constant as desired. This proves Theorem A.

Similarly, Theorem B will be a consequence of the following:

\begin{Theorem} Let $M$ be a closed oriented surface of negative
Gaussian curvature and $\Omega$ and arbitrary 2-form. There exists
$\la_{0}>0$ with the following property. If $G:M\to\re$ is any
smooth function and $\omega$ is any smooth 1-form on $M$ such that
there is a smooth function $g:SM\to\re$ for which
$G(x)+\omega_{x}(v)=X_{\la}(g)$ for $|\la|<\la_{0}$, then $G$ is
constant and $\omega$ is exact. \label{1forma+F}
\end{Theorem}

If we apply Theorem \ref{1forma+F} to (\ref{prekey}) we conclude
as above that $K$ and $F$ must be constant. This proves Theorem B.

\qed

\section{Proof of Theorems \ref{1forma} and \ref{1forma+F}}

Theorems \ref{1forma} and \ref{1forma+F} will be consequences of
more general theorems which are the analogue of Theorem 3.6 in
\cite{GK}. In this section we will try to follow as closely as
possible the notation in \cite{GK}.

Define
\[\eta^{+}:=(X-i\,H)/2\]
and
\[\eta^{-}:=(X+i\,H)/2.\]

Let $L^2(SM)$ be the space of square integrable functions with respect
to the Liouville measure of $SM$.
The next proposition summarizes the main properties of these
operators.

\begin{Proposition}[\cite{GK}] We have:

\begin{enumerate}
\item $L^{2}(SM)$ decomposes into an orthogonal direct sum of
subspaces $\sum H_{n}$, $n\in\Z$, such that on $H_{n}$, $-i\,V$ is
$n$ times the identity operator;
\item $\eta^{+}$ extends to a densely defined operator from
$H_{n}$ to $H_{n+1}$ for all $n$. Moreover, its transpose is
$-\eta_{-}$;
\item let $A:=\min(-K/2)$. Then for all $f\in H_{n}\cap \mbox{\rm
domain}\,\eta^{+}\cap \mbox{\rm domain}\,\eta^{-}$ and $n\geq 0$
\[\|\eta^{+}f\|^{2}\geq An\|f\|^{2}+\|\eta^{-}f\|^{2}.\]

\end{enumerate}
\label{properties}
\end{Proposition}

\begin{Theorem} Let $M$ be a closed oriented surface and let
$N$ be a non-negative integer. Let $\la$ be a real number such
that $\la^2\max\{(N+1),2\}+K(x)<0$ for all $x\in M$. Let
$f:SM\to\re$ be a smooth function of the form
\[f=\sum_{|i|\leq N}f_n,\;\;\;\;\;f_n\in H_n.\]
Suppose the integral of $f$ over every periodic orbit of $\phi^\la_{t}$ is zero.
Then there exists a smooth function $g:SM\to\re$ of the form
\[g=\sum_{|i|\leq N-1}g_n,\;\;\;\;\;g_n\in H_n\]
such that $X_{\la}(g)=f$. (If $N=0$ we interpret this as saying
that $g=0$.) \label{key}
\end{Theorem}

\begin{proof} By the smooth version of the Livsic theorem \cite{LMM}, 
there exists a smooth function $g$
with $X_{\la}(g)=f$. Write
$$g=\sum_{-\infty}^{\infty}g_n,$$
where $g_n\in H_n$. The equation $X_{\la}(g)=f$ is equivalent to the system
of equations
\[\eta^{+}g_{n-1}+\eta^{-}g_{n+1}+in\la g_n=f_n\;\;\;\;\;n=0,\pm 1,\pm 2,\dots.\]
Since $f_n=0$ for $n>N$ we get
\begin{equation}
\eta^{+}g_{n-1}+\eta^{-}g_{n+1}+in\la g_n=0\;\;\;\;\;n>N.
\label{1}
\end{equation}
Using Proposition \ref{properties} item 3 and equation (\ref{1})
we obtain:
\begin{align*}
\|\eta^{+}g_{n+1}\|^2 &\geq \|\eta^{-}g_{n+1}\|^2+A(n+1)\| g_{n+1}\| ^2\\
& =\|\eta^{+}g_{n-1}+in\la g_n\|^2+ A(n+1)\| g_{n+1}\| ^2\\
& =\|\eta^{+}g_{n-1}\|^2+2n\la\mbox{\rm Re}\langle \eta^{+}g_{n-1},i g_n\rangle+
n^2\la^2\| g_n\| ^2+A(n+1)\| g_{n+1}\| ^2.
\end{align*}
For $n>N+1$ we can also write:
\[\|\eta^{+}g_{n}\|^2
\geq \|\eta^{+}g_{n-2}\|^2+2(n-1)\la\mbox{\rm Re}\langle \eta^{+}g_{n-2},i g_{n-1}\rangle+
(n-1)^2\la^2\| g_{n-1}\| ^2+An\| g_{n}\|^2.\]
Thus if we set
\[a_n:= \|\eta^{+}g_{n}\|^2+\|\eta^{+}g_{n-1}\|^2\]
the last two inequalities imply for $n>N+1$:
\begin{align*}
a_{n+1}&\geq a_{n-1}+2n\la\mbox{\rm Re}\langle \eta^{+}g_{n-1},i g_n\rangle+
2(n-1)\la\mbox{\rm Re}\langle \eta^{+}g_{n-2},i g_{n-1}\rangle\\
&+n^2\la^2\| g_n\| ^2+A(n+1)\| g_{n+1}\| ^2+(n-1)^2\la^2\| g_{n-1}\| ^2+An\| g_{n}\|^2.
\end{align*}
Now we compute using equation (\ref{1}) again:
\begin{align*}
\mbox{\rm Re}\langle \eta^{+}g_{n-1},i g_n\rangle&=
\mbox{\rm Re}\langle g_{n-1},-\eta^{-}(i g_n)\rangle\\
&=\mbox{\rm Re}\langle g_{n-1},-i(-i(n-1)\la g_{n-1}-\eta^{+}g_{n-2})\rangle\\
&= -(n-1)\la\| g_{n-1}\|^2 -\mbox{\rm Re}\langle \eta^{+}g_{n-2},i g_{n-1}\rangle
\end{align*}
for $n>N+1$. Hence
\begin{align*}
2n\la\mbox{\rm Re}\langle \eta^{+}g_{n-1},i g_n\rangle+
2(n-1)\la\mbox{\rm Re}\langle \eta^{+}g_{n-2},i g_{n-1}\rangle &\\
=-2n(n-1)\la^{2}\| g_{n-1}\| ^{2}-2\la\mbox{\rm Re}\langle \eta^{+}g_{n-2},i g_{n-1}\rangle.
\end{align*}
This yields for $n>N+1$:
\begin{align*}
a_{n+1}&\geq a_{n-1}-2\la\mbox{\rm Re}\langle \eta^{+}g_{n-2},i g_{n-1}\rangle\\
&+n^2\la^2\| g_n\| ^2+A(n+1)\| g_{n+1}\| ^2-(n-1)(n+1)\la^2\| g_{n-1}\| ^2+An\| g_{n}\|^2.
\end{align*}
For $n>N+2$ we can also write:
\begin{align*}
a_{n}&\geq a_{n-2}-2\la\mbox{\rm Re}\langle \eta^{+}g_{n-3},i g_{n-2}\rangle\\
&+(n-1)^2\la^2\| g_{n-1}\| ^2+An\| g_{n}\|^2-(n-2)n\la^2\| g_{n-2}\| ^2+A(n-1)\| g_{n-1}\|^2.
\end{align*}
Thus if we set
\[b_n:=a_n+a_{n-1}\]
the last two inequalities imply for $n>N+2$:
\begin{align*}
b_{n+1}\geq & b_{n-1}+2(n-2)\la^2\| g_{n-2}\| ^{2}-(n-2)n\la^2\| g_{n-2}\| ^2 \\
&+(n-1)(A-2\la^2)\| g_{n-1}\| ^{2}+(2An+\la^{2}n^{2})\| g_{n}\| ^{2}+A(n+1)\| g_{n+1}\|^2.
\end{align*}
This is the basic recurrence inequality that we will use to prove the theorem.
Let us call
\[r_{n}:=-(n-2)^{2}\la^{2}\| g_{n-2}\| ^{2}
+(n-1)(A-2\la^2)\| g_{n-1}\| ^{2}+(2An+\la^{2}n^{2})\| g_{n}\|
^{2}\]
\[+A(n+1)\| g_{n+1}\|^2\] which is defined for $n>N+2$. We have:
\begin{equation}
b_{n+1}\geq b_{n-1}+r_{n}. \label{recu}
\end{equation}
Consider now an positive integer $n$ of the form $n=N+3+2k$ for
$k\geq 0$ and let $j$ be another integer with $k\geq j\geq 0$.
Using (\ref{recu}) we obtain:
\begin{equation}b_{n+1}\geq
b_{N+2+2j}+r_n+r_{n-2}+\dots+r_{N+3+2j}.\label{recu1}
\end{equation} Note that
\[r_{n}+r_{n-2}=-(n-4)^{2}\la^{2}\| g_{n-4}\|
^{2}+(n-3)(A-2\la^{2})\| g_{n-3}\| ^{2}+2A(n-2)\| g_{n-2}\| ^{2}\]
\[+2(n-1)(A-\la^2)\| g_{n-1}\| ^{2}+(2An+\la^{2}n^{2})\| g_{n}\|
^{2}+A(n+1)\| g_{n+1}\|^2.\] If we now suppose that $A-\la^2\geq
0$, then a simple calculation gives:
\[r_n+r_{n-2}+\dots+r_{N+3+2j}\geq -(N+1+2j)^2\la^{2}\| g_{N+1+2j}\|^2\]
\[+(N+2+2j)(A-2\la^2)
\| g_{N+2+2j}\|^2.\] Using (\ref{recu1}) we obtain:
\[b_{n+1}\geq \| \eta^{+}g_{N+2+2j}\| ^{2}+2\| \eta^{+}g_{N+1+2j}\| ^{2}
+\| \eta^{+}g_{N+2j}\| ^{2}\] \[-(N+1+2j)^2\la^{2}\|
g_{N+1+2j}\|^2+(N+2+2j)(A-2\la^2)\| g_{N+2+2j}\|^2,\] which
combined with Proposition \ref{properties} item 3 gives:
\[b_{n+1}\geq 2(N+2+2j)(A-\la^2)\| g_{N+2+2j}\| ^{2}\]
\begin{equation}+(N+1+2j)(2A-(N+1+2j)\la^{2})\| g_{N+1+2j}\| ^{2}+ A(N+2j)\|
g_{N+2j}\| ^{2}. \label{2}
\end{equation}
Analogously, if we take a positive integer $n$ of the form
$n=N+3+2k+1$ for $k\geq 0$ and let $j$ be another integer with
$k\geq j\geq 0$ we obtain:
\[b_{n+1}\geq 2(N+3+2j)(A-\la^2)\| g_{N+3+2j}\| ^{2}\]
\begin{equation}+(N+2+2j)(2A-(N+2+2j)\la^{2})\| g_{N+2+2j}\| ^{2}+ A(N+1+2j)\|
g_{N+1+2j}\| ^{2}. \label{3}
\end{equation}

Suppose now that $A-\la^2>0$ and $2A-(N+1)\la^{2}>0$. Since $g$ is
a smooth function, the functions $g_n$ must tend to zero in the
$L_{2}$-topology as $n$ tends to infinity. Thus $b_n\to 0$, which
can only occur if $g_{N}=g_{N+1}=g_{N+2}=0$ in view of (\ref{2})
for $j=0$. If we now use that $g_{N+1}=g_{N+2}=0$ in (\ref{3}) for
$j=0$, the fact that $b_{n}\to 0$ shows that $g_{N+3}=0$. If we
continue in this fashion using (\ref{2}) and (\ref{3}) for all
values of $j$ we get that $g_{n}=0$ for all $n\geq N$.

Analogously one can prove that $g_n=0$ for $n\leq -N$.

\end{proof}

Theorem \ref{1forma} follows right away from Theorem \ref{key}:
given a 1-form $\omega$ on $M$, regarded as a function
$\omega:TM\to\re$, we can consider its restriction to $SM$ which
lies in $H_{1}\oplus H_{-1}$. (If we write
$\omega=\omega_{1}+\omega_{-1}$, then
$\bar{\omega}_{1}=\omega_{-1}$.) Theorem \ref{key} tell us that
there exists $g\in H_0$ such that $\omega=X_{\la}(g)$. But if
$g\in H_0$, then it is the pull back of a smooth function
$h:M\to\re$ and thus $X_{\la}(g)=dh$ and $\omega$ is exact as
desired.

Similarly, Theorem \ref{1forma+F} follows immediately from the
following theorem. Its proof is quite similar to that of Theorem
\ref{key}, so we will only indicate the main steps.

\begin{Theorem} Let $M$ be a closed oriented surface of negative Gaussian curvature
and let $\Omega$ be an arbitrary smooth 2-form. Let $N$ be a
non-negative integer. There exists a positive number $\lambda_{0}$
which depends on the metric, $\Omega$ and $N$ such that for all
$\la\in [0,\la_{0})$ the following property holds: let
$f:SM\to\re$ be a smooth function of the form
\[f=\sum_{|i|\leq N}f_n,\;\;\;\;\;f_n\in H_n.\]
Suppose the integral of $f$ over every periodic orbit of
$\phi^\la_{t}$ is zero. Then there exists a smooth function
$g:SM\to\re$ of the form
\[g=\sum_{|i|\leq N-1}g_n,\;\;\;\;\;g_n\in H_n\]
such that $X_{\la}(g)=f$. (If $N=0$ we interpret this as saying
that $g=0$.) \label{keygen}
\end{Theorem}

\begin{proof} We will suppose without loss of generality that the metric is
normalized so that $K\leq -2$.

Using
\[\eta^{+}g_{n-1}+\eta^{-}g_{n+1}+in\la\,F\, g_n=0,\;\;\;\;\;n>N, \]
the same proof of Theorem \ref{key} shows that (with the same
definition of $b_n$) for $n>N+2$ we have:

\begin{align*}
b_{n+1}\geq & b_{n-1}-2(n-1)\|\la\,F\,
g_{n-1}\|^2-(n-2)^{2}\|\la\,F\, g_{n-2}\| ^{2}+n^{2}\|\la\,F\,
g_{n}\| ^2 \\ &+(n-1)A\| g_{n-1}\| ^{2}+2An\| g_{n}\|
^{2}+A(n+1)\| g_{n+1}\|^2\\ & +2(n-2)\,\mbox{\rm Re}\langle
g_{n-2},-i\eta^{-}(\la\,F)g_{n-1}\rangle+2n\,\mbox{\rm Re}\langle
g_{n-1},-i\eta^{-}(\la\,F)g_{n}\rangle.\\
\end{align*}
We now observe:
\[2(n-2)\,\mbox{\rm Re}\langle
g_{n-2},-i\eta^{-}(\la\,F)g_{n-1}\rangle+2n\,\mbox{\rm Re}\langle
g_{n-1},-i\eta^{-}(\la\,F)g_{n}\rangle\]
\[\geq -(n-2)\| g_{n-2}\|
^{2}-2(n-1)\|\eta^{-}(\la\,F) g_{n-1}\| ^{2}-n\| g_{n}\| ^{2}.\]
Since $M$ is compact there exists a positive constant $c$ such
that for all $n$ and $g$ we have
\[\|\eta^{-}(F) g_{n}\| ^{2}\leq c\,\| g_{n}\| ^{2},\]
\[\|F\, g_{n}\| ^{2}\leq c\,\| g_{n}\| ^{2}\]
and thus

\begin{align*}
b_{n+1}\geq & b_{n-1}-(n-2)^{2}\|\la\,F\, g_{n-2}\|
^{2}+n^{2}\|\la\,F\, g_{n}\| ^2-(n-2)\| g_{n-2}\| ^{2} \\
&+(n-1)(A-4c\la^{2})\| g_{n-1}\| ^{2}+(2A-1)n\| g_{n}\|
^{2}+A(n+1)\| g_{n+1}\|^2.\\
\end{align*}
By choosing $\la$ small enough we may suppose that
$A-4c\la^{2}\geq 0$. Our normalization of the metric says that
$A\geq 1$ and hence

\[b_{n+1}\geq  b_{n-1}-(n-2)^{2}\|\la\,F\, g_{n-2}\|
^{2}+n^{2}\|\la\,F\, g_{n}\| ^2-(n-2)\| g_{n-2}\| ^{2}+n\| g_{n}\|
^{2}.\] This implies
\[b_{n+1}\geq b_{N+2}-(N+1)^{2}\|\la\,F\, g_{N+1}\| ^2-(N+1)\| g_{N+1}\|
^{2}.\] Using the definition of $b_n$ and Proposition
\ref{properties} item 3 we obtain:
\[b_{n+1}\geq A(N+2)\| g_{N+2}\|
^{2}+[(2A-1)(N+1)-\la^{2}c(N+1)^{2}]\| g_{N+1}\| ^{2}+AN\| g_{N}\|
^{2}\] and if we now choose $\la$ small enough so that
\[(2A-1)(N+1)-\la^{2}c(N+1)^{2}>0\]
then the same argument as in Theorem \ref{key} shows that
$g_N=g_{N+1}=g_{N+2}=0$. Now it follows easily that $g_n=0$ for
all $n\geq N$.

\end{proof}

\begin{Remark}{\rm It is quite likely that one can prove (and perhaps improve)
Theorems \ref{key} and \ref{keygen} by proving first a Pestov's identity
for magnetic flows and then proceeding as in \cite{CS,DS}.}
\end{Remark}

\end{document}